\documentclass[a4paper,12pt]{article}
\usepackage{graphicx}
\usepackage{a4wide}
\usepackage[english]{babel}
\usepackage{amsfonts}
\usepackage{amsmath}
\usepackage{amssymb}
\usepackage{dsfont}

\newtheorem{lemma}{Lemma}

\newtheorem{theorem}{Theorem}
\newtheorem{corollary}{Corollary}

\allowdisplaybreaks

\date{}

 \DeclareMathOperator
{\sgn}{sgn}

\newcommand {\e} {\varepsilon}

\newcommand {\Z} {\mathbb{Z}}

\newcommand {\ve} {\varepsilon}

\numberwithin{equation}{section} \numberwithin{theorem}{section}
\numberwithin{lemma}{section} 
\numberwithin{corollary}{section}
\makeatletter\def\blfootnote{\xdef\@thefnmark{}\@footnotetext}
\makeatother

\begin{document}
\title{\bf On the system $f(nx)$ and probabilistic number theory}
\author{C.\ Aistleitner\footnote{Graz University of Technology,
Institute of Mathematics A, Steyrergasse 30, 8010 Graz, Austria.
\mbox{e-mail}: \texttt{aistleitner@math.tugraz.at}. Research
supported by FWF grant S9610-N23.}, I.\ Berkes\footnote{ Graz
University of Technology, Institute of Statistics,
M\"unzgrabenstra{\ss}e 11, 8010 Graz, Austria. \mbox{e-mail}:
\texttt{berkes@tugraz.at}. Research supported by the FWF Doctoral Program on Discrete Mathematics (FWF DK W1230-N13), FWF grant
S9603-N23 and OTKA grant K 81928.} and R.\ Tichy\footnote{Graz
University of Technology, Institute of Mathematics A, Steyrergasse
30, 8010 Graz, Austria. \mbox{e-mail}: \texttt{tichy@tugraz.at}.
Research supported by the FWF Doctoral Program on Discrete Mathematics (FWF DK W1230-N13) and FWF grant S9603-N23.}}

\maketitle \vskip1cm

\begin{center}
\emph{Dedicated to the memory of Professor Jonas Kubilius.}
\end{center}

\vskip1cm

\abstract{Let $f: {\mathbb R}\to {\mathbb R}$ be a measurable
function satisfying
\begin{equation*}
f(x+1)=f(x), \qquad   \int_0^1 f(x)\, dx=0, \qquad \int_0^1
f^2(x)\, dx<\infty.
\end{equation*}
The asymptotic properties  of series  $\sum c_k f(kx)$ have been
studied extensively in the literature and turned out to be, in
general, quite different from those of the trigonometric system.
As the theory shows, the behavior of such series is determined by
a combination of analytic, probabilistic and number theoretic
effects, resulting in highly interesting phenomena not encountered
in classical harmonic analysis. In this paper we survey some
recent results in the field and prove asymptotic results for the
system $\{f(nx), n\ge 1\}$ in the case when the function $f$ is
not square integrable.



\newpage
\section{Introduction}

Let $f: {\mathbb R}\to {\mathbb R}$ be a measurable function
satisfying
\begin{equation}\label{f}
f(x+1)=f(x), \qquad   \int_0^1 f(x)\, dx=0, \qquad \int_0^1
f^2(x)\, dx<\infty.
\end{equation}
The asymptotic properties  of the system $\{f(nx), n\ge 1\}$ have
been studied extensively in the literature and turned out to be,
in general, very different from those of the trigonometric system.
Khinchin \cite{khin2} conjectured that (\ref{f}) (even without the
last condition) implies
\begin{equation}\label{h}
\lim_{N\to \infty} \frac{1}{N} \sum_{k=1}^N f(kx)=0 \qquad
\hbox{a.e.}
\end{equation} This remained open for almost 50 years
until it was disproved by Marstrand \cite{mar}; he showed that
there exist even bounded counterexamples $f$. At about the same
time, Nikishin \cite{nik} constructed a continuous function $f$
satisfying (\ref{f}) such that $\sum c_k f(kx)$ diverges a.e.\ for
some $(c_k)$ with $\sum c_k^2<\infty$. Gaposhkin \cite{gapo1968b}
showed that if $f$ satisfies (\ref{f}) and belongs to the
Lipschitz $\alpha$ class for some $\alpha>1/2$, then $\sum c_k
f(kx)$ converges a.e.\ provided $\sum c_k^2<\infty$, i.e. the
analogue of Carleson's theorem holds for the system $f(nx)$.
Berkes \cite{b98} showed that this result becomes false, in
general, for Lipschitz 1/2 functions. There exists no
characterization of functions $f$ for which the analogue of the
Carleson convergence theorem holds for $f(nx)$ and, despite the
profound work of Bourgain \cite{bou} connecting Khinchin's
conjecture with metric entropy behavior, we have no
characterization of functions $f$ for which (\ref{h}) holds.

The asymptotic properties of the system $\{f(nx), n\ge 1\}$ play
also an important role in the metric theory of uniform
distribution. A sequence $(x_n)_{n\ge 1}$ of real numbers is
called uniformly distributed mod 1 if for any interval $[a, b)$ on
the real line we have
$$\lim_{N\to \infty} \frac{1}{N} \sum_{k=1}^N
\mathds{1}_{[a,b)} (x_k)=b-a,
$$
where  $\mathds{1}_{[a,b)}$  denotes the indicator function of the
interval $[a,b)$, extended with period 1. Given a sequence $(x_1,
\dots, x_N)$ of real numbers, the value
$$
D_N = D_N (x_1, \dots, x_N) = \sup_{0 \leq a < b < 1} \left|
\frac{\sum_{k=1}^N \mathds{1}_{[a,b)} (x_k)}{N} - (b-a) \right|
$$
is called the discrepancy of the sequence. It is easily seen that
an infinite sequence $(x_n)_{n\ge 1}$ is uniformly distributed mod
1 iff $D_N(x_1, \ldots x_N) \to 0$ as $N\to \infty$. By a
classical result of Weyl \cite{we16}, for any increasing sequence
$(n_k)$ of integers, $\{n_kx\}_{k\ge 1}$ is uniformly distributed
mod 1 for all $x\in{\mathbb R}$, with the exception of a set
having Lebesgue measure 0. Improving the results of Erd\H{o}s and
Koksma \cite{erko} and Cassels \cite{ca}, Baker \cite{bak} proved
that for any increasing sequence $(n_k)$ of positive integers the
discrepancy $D_N (\{n_kx\})$  of the first $N$ terms of the
sequence $\{n_kx\}$ satisfies
\begin{equation}\label{baker}
D_N (\{n_kx\}) \ll \frac{ (\log N)^{3/2 +\varepsilon}}{\sqrt{N}}
\qquad \text{a.e.}
\end{equation}
for any $\varepsilon>0$. On the other hand, Berkes and Philipp
\cite{bp} constructed an increasing sequence $(n_k)$ of integers
such that for almost all $x$ the relation
\begin{equation}\label{bp}
D_N (\{n_kx\}) \ge \frac{ (\log N)^{1/2}}{\sqrt{N}}
\end{equation}
holds for infinitely many $N$. These results describe the extremal
behavior of $D_N(\{n_kx\})$ rather precisely; on the other hand,
the exact asymptotics of $D_N(\{n_kx\})$ is known only in a few
special cases, e.g.\ for $n_k=k$ (Khinchin \cite{kh24}, Kesten
\cite{ke}) and for exponentially growing $n_k$ (Philipp
\cite{plt}). Let us note that for every $x\in {\mathbb R}$
\begin{equation}\label{est}
\frac{1}{4} \sup_{V_f\le 2}\left| \sum_{k=1}^N f(n_kx)\right| \le
D_N(\{n_kx\}) \le\sup_{V_f\le 2} \left| \sum_{k=1}^N
f(n_kx)\right|,
\end{equation}
where $V_f$ denotes the total variation of $f$ on $[0, 1]$. The
second inequality in (\ref{est}) is obvious from the definition of
$D_N(\{n_kx\})$, while the first one  follows from Koksma's
inequality (see e.g. \cite{dts}), stating
\begin{equation} \label{koksma}
\left| \frac{1}{N} \sum_{k=1}^N f(x_k) - \int_0^1 f(x)~dx \right|
\leq 2 V_f \cdot D_N (x_1, \ldots x_N)
\end{equation}
for any function $f$ with $V_f<\infty$  and for every set $x_1,
\dots, x_N$ of points from the unit interval. The inequality
(\ref{koksma}) plays a crucial role in the theory of Monte Carlo
and quasi-Monte Carlo integration. By (\ref{est}), determining the
precise order of magnitude of $D_N(\{n_kx\})$ requires sharp
bounds for sums $\sum_{k=1}^N f(n_kx)$, which, in turn, is closely
connected with estimating the integral
\begin{equation}\label{int}
 \int_0^1 \left(\sum_{k=1}^N
f(n_kx)\right)^2 dx
\end{equation}
for functions $f$ of bounded variation. Koksma \cite{ko51} showed
that the integral in (\ref{int}) is bounded by $V_f \cdot G(n_1, n_2,
\ldots, n_N)$, where
\begin{equation}
G(n_1, n_2, \ldots, n_N)= \sum_{1\le i\le j \le N} \frac{ (n_i,
n_j)}{[n_i, n_j]}.
\end{equation}
This shows  that the order of magnitude of the discrepancy
$D_N(\{n_kx\})$ depends not only on the growth speed of the
sequence, but also its number theoretic properties.

Relations (\ref{baker}) and (\ref{est}) imply that for any $f$
with $V_f<\infty$ and any increasing sequence $(n_k)$ of integers
we have
\begin{equation}\label{fest}
\left| \sum_{k=1}^N f(n_kx)\right| \ll \sqrt{N} (\log N)^{3/2
+\ve} \qquad \text{a.e.}
\end{equation}
for any $\ve>0$. On the other hand, the proof of (\ref{bp}) in
Berkes and Philipp \cite{bp} provides an increasing sequence
$(n_k)$ of integers such that for $f(x)=\{x\}-1/2$ we have for
almost all $x$
\begin{equation}\label{bp5}
\left| \sum_{k=1}^N f(n_kx)\right| \ge \sqrt{N} (\log N)^{1/2}
\end{equation}
for infinitely many $N$. The gap between (\ref{baker}) and
(\ref{bp}), as well as the gap between (\ref{fest}) and
(\ref{bp5}) remain open until today. Very recently, Aistleitner,
Mayer and Ziegler \cite{aizi} improved (\ref{fest}) to
\begin{equation}\label{az}
\left|\sum_{k=1}^N f(n_kx)\right|\ll \sqrt{N} (\log N)^{3/2}
(\log\log N)^{-1/2+\varepsilon} \qquad \text{a.e.}
\end{equation}
for any $\ve >0$. Their proof uses the estimate
\begin{equation}\label{har}
\sum_{1\le k\le l\le N} \frac{(n_k, n_l)}{\sqrt{n_k n_l}} \ll N
\exp \left(\frac{c\log N}{\log\log N}\right)
\end{equation}
due to Dyer and Harman \cite{dyha}. Dyer and Harman also
conjectured that $\log N$ on the right hand side of (\ref{har})
can be replaced by $\sqrt{\log N}$. Assuming the validity of this
conjecture, the bound in (\ref{az}) can be improved to
$\sqrt{N}(\log N)^{1+\e}$.

\medskip
The previous results show that there is a crucial difference
between the asymptotic behavior of the system $f(nx)$ for
functions $f$ with bounded variation and for $f(x)=e^{2\pi ix}$.
In the latter case Berkes and Philipp \cite{bp} proved that if
$\psi (n)$ is a nondecreasing sequence satisfying $\psi(n^2)\ll
\psi(n)$, then
\begin{equation}\label{bp3} \left|\sum_{k=1}^N e^{2\pi i
n_kx}\right| \ll \sqrt{N} \psi (N) \qquad \text{a.e.}
\end{equation}
holds for all increasing $(n_k)$ if and only if
$$\sum_{k=1}^\infty \frac{1}{k\psi (k)^2}<\infty.$$
In particular, (\ref{bp3}) holds for $\psi (N)=(\log N)^{1/2+\ve}$
for $\ve>0$ and fails for  $\psi (N)=(\log N)^{1/2}$

\medskip
By a profound result of G\'al \cite{gal}, for any sequence $(n_1,
\ldots, n_N)$ we have
\begin{equation}\label{gal}
G(n_1, n_2, \ldots, n_N) \ll N(\log \log N)^2
\end{equation}
and this result is best possible. Since under $V_f<\infty$ the
integral (\ref{int}) is $\ll G(n_1, n_2, \ldots, n_N)$, for
functions $f$ with bounded variation the $L^2$ norm of $\sum_{k=1}^N
f(n_kx)$ is $O(\sqrt{N} \log\log N)$, a bound only slightly weaker
than the bound $O(\sqrt{N})$ valid for orthogonal series. Thus one
can expect that under $V_f<\infty$ the convergence properties of
$\sum_{k=1}^\infty  c_k f(n_kx)$ are also not much worse than those
of orthogonal series, described by the Rademacher-Mensov convergence
theorem.  This is indeed the case: Berkes and Weber \cite{bewe09}
showed that for any increasing sequence $(n_k)$ of positive integers
and any function $f$ satisfying (\ref{f}) and $V_f<\infty$, the
series $\sum_{k=1}^\infty c_k f(n_kx)$ converges a.e.\ provided
\begin{equation}\label{rm}
\sum_{k=1}^\infty  c_k^2 (\log k)^{3+\ve}  <\infty
\end{equation}
for some $\ve>0$. On the other hand, Nikishin \cite{nik} showed that
for $f(x)=\sgn \sin 2\pi x$ (a function with bounded variation on
$[0, 1]$) the series $\sum_{k=1}^\infty c_k f(kx)$ diverges on a set
with positive measure for some $(c_k)$ with $\sum_{k=1}^\infty
c_k^2<\infty$. These two results characterize, up to a logarithmic
factor, the a.e.\ convergence of $\sum_{k=1}^\infty c_k f(n_kx)$ for
functions $f$ with bounded variation. For other classes of functions
the convergence properties of the series are completely different.
Recall that if $f$ is a Lip $\alpha$ function with $\alpha>1/2$
satisfying (\ref{f}), then the analogue of Carleson's theorem holds
for $\sum_{k=1}^\infty c_k f(kx)$, and this theorem generally fails
for $\alpha=1/2$. In the case $0<\alpha<1/2$ Weber \cite{we11}
proved, improving results of Gaposhkin \cite{gapo1967}, that a
sufficient convergence criterion is
\begin{equation}\label{wd}
\sum_{k=1}^\infty  c_k^2 d(k) (\log k)^2<\infty,
\end{equation}
where $d(k)=\#\{1\le i\le k: i|k\}$  is the divisor function. It
is known that
\begin{equation}\label{dk} d(k)\ll \exp (C\log k/\log \log
k),
\end{equation}
for some $C>0$ and thus $\sum_{k=1}^\infty c_k f(kx)$  converges
a.e.\ provided
$$\sum c_k^2 \exp (C_1\log k/\log\log k)<\infty,$$
a fact proved independently also by Aistleitner \cite{ai11} for
$1/4<\alpha<1/2$. We note that, as Weber \cite{we11} showed,
(\ref{wd}) is sufficient for the a.e.\ convergence of
$\sum_{k=1}^\infty c_k f(kx)$ even if instead of the Lipschitz
character of $f$ we assume only that the Fourier coefficients of
$f$ are $O(k^{-1/2} (\log k)^{-(1+\ve)})$ for some $\ve >0$, a
criterion allowing a much larger class of functions $f$. In the
case when the Fourier coefficients of $f$ are $O(k^{-\gamma})$,
$1/2<\gamma<1$, a sufficient convergence criterion is
$$\sum_{k=1}^\infty c_k^2 \rho_\gamma (k) (\log k)^2 <\infty$$
where
$$\rho_\gamma (n)=\sum_{d|n} d^{-(2\gamma-1)}.
$$
See Berkes and Weber \cite{bewe12}. These results show that even
in the case $n_k=k$ the convergence behavior of $\sum_{k=1}^\infty
c_k f(n_kx)$ is intimately connected with number theory. For a
detailed study of the convergence properties of sums
$\sum_{k=1}^\infty c_k f(n_kx)$, see Berkes and Weber
\cite{bewe09}.

The previous results show that even for ``nice'' functions $f$, the
a.e.\ convergence of $\sum_{k=1}^\infty c_k f(kx)$ is a highly
delicate question, far from being solved. In contrast, convergence
in $L^2$ norm is essentially solved by a theorem of Wintner
\cite{win}, who proved that if $f$ has the Fourier series
$$f\sim \sum_{k=1}^\infty  (a_k \cos 2\pi kx+ b_k \sin 2\pi kx),$$
then  $\sum_{k=1}^\infty c_k f(kx)$ converges in norm for all
$(c_k)$ with $\sum_{k=1}^\infty  c_k^2<\infty$ iff the Dirichlet
series
$$
\sum_{k=1}^\infty a_k k^{-s}, \quad \hbox{and} \quad
\sum_{k=1}^\infty b_k k^{-s} \eqno(1.4)
$$
are regular and bounded in the half-plane ${\Re (s)}>0$. This
remarkable criterion shows again the complexity of the convergence

problem studied here.

\section{Lacunary series}

In the previous chapter we have seen that that the convergence and
growth properties of series $\sum_{k=1}^\infty  c_k f(n_kx)$ are
closely connected with the number theoretic properties of $(n_k)$,
and even for $n_k=k$ the convergence problem has an arithmetic
character. In this chapter we investigate lacunary series, i.e.
the behavior of $f(n_kx)$, where $(n_k)$ satisfies the Hadamard
gap condition
\begin{equation}\label{had}
n_{k+1}/n_k \ge q >1 \qquad (k=1, 2, \ldots).
\end{equation}
In the case $f(x)=\cos 2\pi x$ and $f(x)=\sin 2\pi x$, Salem and
Zygmund \cite{salz} and Erd\H{o}s and G\'al \cite{erga} proved the
central limit theorem  and the law of the iterated logarithm,
i.e.\
\begin{equation}\label{clt}
\frac{1}{\sqrt{N/2}}\sum_{k=1}^N
f(n_kx)\buildrel{d}\over{\longrightarrow} N(0, 1)
\end{equation}
and
\begin{equation}\label{lil}
\limsup_{N\to \infty} \frac{1}{\sqrt{N\log \log N}}\sum_{k=1}^N
f(n_kx)=1 \qquad \text{a.e.}
\end{equation}
with respect to the probability space $([0, 1], {\mathcal B},
\mu)$,  where $\mathcal B$ is the Borel $\sigma$-field in $[0, 1]$
and $\mu$ is the Lebesgue measure. This shows that Hadamard
lacunary subsequences of the trigonometric system behave like
independent random variables. Extensions for weighted sums
$\sum_{k=1}^N c_k f(n_kx)$ were proved by Salem and Zygmund
\cite{salz,salz2} and  Weiss \cite{weiss} under the same
coefficient conditions as assumed  for independent random
variables. For general $f$ satisfying (\ref{f}), the situation is
considerably more complex. Kac \cite{kac1946} showed that if $f$
is a Lipschitz function or a function with bounded variation
satisfying (\ref{f}), then the CLT (\ref{clt}) holds in the case
$n_k=2^k$ with a limit distribution $N(0, \sigma^2)$, where
\begin{equation}\label{var}
\sigma^2=\int_0^1 f^2(x)\, dx+2\sum_{k=1}^\infty \int_0^1
f(x)f(2^kx) \, dx.
\end{equation}
The corresponding law of the iterated logarithm was proved by
Izumi \cite{iz} and Maruyama \cite{maru}. On the other hand, Erd\H
os and Fortet showed (see \cite{kac1949}, p.\ 646) that both the
CLT and LIL fail for the system $f(n_k x)$, where
$$
f(x)=\cos 2 \pi x + \cos 4 \pi x \quad \textrm{and} \quad n_k =
2^k -1, \quad k \geq 1.
$$
Specifically, in the case of the CLT a limit distribution in
(\ref{clt}) still exists, but its distribution function  equals
$$\pi^{-1/2} \int_0^1 \int_{-\infty}^{x/2|\cos \pi t|}e^{-u^2}dudt$$
(i.e.\ it is a mixture of Gaussian distributions) and the limsup
in  (\ref{lil}) is $\sqrt{2} \cos \pi x$, i.e.\ the limsup depends
on $x$ (cf. also \cite{conze}). These results show that under (\ref{had}) the behavior of
$f(n_kx)$ still resembles to that of independent random variables,
but it is influenced substantially by the number theoretic
properties of the sequence $(n_k)$ as well. Gaposhkin
\cite{gapo1966} showed that the CLT for $f(n_kx)$ remains valid if
all the fractions $n_{k+1}/n_k$ are integers or if $n_{k+1}/n_k
\to \alpha$, where $\alpha^r$ is irrational for $r=1, 2, \ldots$.
He also showed (see \cite{gapo1970}) that the validity of the CLT
is intimately connected with the number of solutions $(k,l)$ of
Diophantine equations of the form
\begin{equation} \label{diophequ}
a n_k \pm b n_l = c, \qquad \textrm{where} \ a,b,c \in \Z .
\end{equation}
Improving Gaposhkin's results, Aistleitner and Berkes \cite{aibe}
gave a complete characterization for the CLT under (\ref{had}).
They proved, namely,  that under (\ref{had}) $f(n_kx)$ satisfies
the CLT for all functions $f$ satisfying (\ref{f}) if and only if
\begin{equation}\label{dlil}
L(N,d,\nu)= o(N) \qquad \textrm{as} \qquad N \to \infty
\end{equation}
uniformly in $\nu\ne 0$, where
\begin{equation}\label{Ldef}
L(N,d,\nu) = \# \{1 \leq a,b \leq d,~1 \leq k,l \leq N:~a n_k - b
n_l = \nu\},
\end{equation}
where we exclude the trivial solutions $k=l$ in the case
$a=b,~\nu=0$. Allowing also $\nu=0$ in (\ref{dlil}), the CLT will
hold with norming  factor $\|f\| \sqrt{N}$. A similar criterion
holds for the LIL, see \cite{ai,aita}.

In his classical paper, Philipp \cite{plt} investigated the law of
the iterated logarithm for the discrepancy $D_N (\{n_kx\})$ under
the lacunarity condition (\ref{had}).   For an i.i.d.\ sequence
$(\xi_n)$ of random variables, uniformly distributed on $(0, 1)$,
the Chung-Smirnov law of the iterated logarithm (see e.g.
\cite{showe}, p.\ 504) states
\begin{equation}\label{chs}
\limsup_{n\to\infty} \frac{N D_N(\xi_1, \ldots \xi_N)}
{\sqrt{2N\log\log N}}=\frac{1}{2}
\end{equation}
with probability 1. Philipp proved that under (\ref{had}) we have
\begin{equation} \label{philipp}
\frac{1}{4 \sqrt{2}} \leq \limsup_{N \to \infty} \frac{N D_N(
\{n_k x\})}{\sqrt{2 N \log \log N}} \leq C_q \qquad \textup{a.e.}
\end{equation}
where the constant $C_q$ depends only on the growth factor $q$ in (\ref{had}). Again,  this shows
that under (\ref{had}) the sequence $\{n_k x\}$ behaves like a
sequence of independent random variables, but the value of the
limsup in (\ref{philipp}), as well as the question whether the
limsup is a constant almost everywhere, remained open. Berkes and
Philipp \cite{bp} showed that for any $q>1$ there exists a
sequence $(n_k)$ satisfying (\ref{had}) such that the limsup in
(\ref{philipp}) exceeds $c \log\log \frac{1}{q}$, showing that the
limsup in (\ref{philipp}) can be different from the classical
value $1/2$ in (\ref{chs}). Aistleitner proved that the $\limsup$
does not have to be a constant a.e.\ for lacunary $(n_k)$ (cf.
\cite{airr1,airr2,aclass}). The (nonconstant)
limsup in the case $n_k=2^k-1$ was determined by Fukuyama
\cite{fna2}.

Very recently Fukuyama \cite{ft} developed a powerful technique to
calculate the exact value  of the limsup in (\ref{philipp}). In
particular, he showed
$$
\limsup_{N \to \infty} \frac{N D_N( \{\theta^k x \})} {\sqrt{2 N
\log \log N}} = C_\theta \quad \textup{a.e.},
$$
where the constants $C_\theta$ (which are explicitly known) depend
on the precise value and the number-theoretic properties of
$\theta$ in a highly interesting way. For example, we have
\begin{eqnarray} \label{Fukuyama}
C_\theta & = & \sqrt{42}/9 , \quad \textrm{if} ~ \theta =2
\nonumber\\
C_\theta & = & \frac{\sqrt{(\theta+1)\theta(\theta-2)}}{2
\sqrt{(\theta-1)^3}} \quad \textrm{if} ~ \theta \geq 4 ~
\textrm{is an even integer},\nonumber\\
C_\theta & = & \frac{\sqrt{\theta+1}}{2 \sqrt{\theta-1}} \quad
\textrm{if} ~ \theta \geq 3 ~ \textrm{is an odd integer.}
\nonumber
\end{eqnarray}
Of particular interest is the case when $\theta$ has no rational
powers (which is the case e.g.\ if $\theta$ is transcendental),
where we have $C_\theta=1/2$, i.e.\ the same constant as in
(\ref{chs}). Aistleitner \cite{ai} showed that the limsup is also
equal $1/2$ if the counting function defined in (\ref{Ldef})
satisfies
\begin{equation}\label{dlil2}
L(N,d,\nu)= O(N/(\log N)^{1+\varepsilon}) \qquad \textrm{as}
\qquad N \to \infty
\end{equation}
for some $\varepsilon>0$, uniformly in $\nu \in {\mathbb Z}$. That
is, under a condition only slightly stronger than the necessary
and sufficient condition for the CLT for $f(n_kx)$,  the
discrepancy behavior of $\{n_kx\}$ also follows i.i.d.\ behavior
precisely. Note, however, that although the asymptotic order of the discrepancy of $\{n_k x\}$ is very well understood
for exponentially growing $(n_k)$ from a probabilistic point of view,
there exist hardly any results for the corresponding problem
for concrete values of $x$. For example, the classical problems asking for uniform distribution of the sequences
$\{2^k \sqrt{2}\}$ and $\{(3/2)^k\}$ are still completely open, and there is little
hope that they can be solved within the next decades (for a discussion of these problems
and recent contributions, see Bailey and Crandall \cite{bc} and Dubickas \cite{dub1,dub2}).\\

Let $\xi_1, \xi_2, \ldots$ be i.i.d.\ random variables, uniformly
distributed on $[0, 1]$, let
$$F_n(x)=\frac{1}{n} \sum_{k=1}^n \mathds{1} (\xi_k \leq x)$$
denote the empirical distribution function of the sample $(\xi_1,
\ldots, \xi_n)$ and let
$$ T_n = \sup_{0\le x \le 1} \sqrt{n} |F_n(x)-x|$$
be the Kolmogorov-Smirnov statistic. $T_n$ plays an important role
in nonparametric statistics, see e.g.\ \cite{showe}. Using
probabilistic terminology, the Chung-Smirnov LIL (\ref{chs}) can
be formulated as
\begin{equation}\label{chs2}
\limsup_{n\to\infty} \frac{T_n}{\sqrt{2\log\log n}}=\frac{1}{2}
\qquad \text{a.s.}
\end{equation}
The limit distributional behavior of $T_n$ is described by
Kolmogorov's  theorem
\begin{equation}\label{k1}
\lim_{n \to \infty} P(T_n\le t)=K(t),
\end{equation}
where 
\begin{equation}\label{k2}
K(t)= 1-2\sum_{k=1}^\infty (-1)^{k-1} e^{-2k^2t^2}.
\end{equation}
It is natural to ask if an analogue of the Kolmogorov limit
theorem (\ref{k1})-(\ref{k2}) holds for discrepancies. Aistleitner
and Berkes \cite{aibe11} showed that if (\ref{dlil}) holds
uniformly in $\nu \in {\mathbb Z}$ (including $\nu=0$), then the
limit distribution of $\sqrt {N} D_N(\{n_kx\})$ exists, namely we
have
\begin{equation*}
\sqrt{N} D_N (\{n_ky\}) \buildrel{d}\over{\longrightarrow} K.
\end{equation*}
Note that the just mentioned Diophantine condition is not
satisfied for $n_k=a^k$ (in this case relation (\ref{dlil}) fails
for $\nu=0$), but the limit distribution of $\sqrt {N}
D_N(\{a^kx\})$ still exists; the limit distribution is the same as
the distribution of $\sup_{0\le t \le 1} |G(x)|$, where $G$ is a
Gaussian process with covariance function
\begin{equation}\label{Gamma}
\Gamma (s, t)=\int_0^1 {\bf I}_s(x) {\bf I}_t(x)\, dx +
\sum_{k=1}^\infty \int_0^1 \left( {\bf I}_s(x) {\bf I}_t (a^k x) +
{\bf I}_s(a^k x) {\bf I}_t (x) \right)\,dx
\end{equation}
where
$$
{\bf I}_t (x) = \mathds{1}_{[0,t]} (x) - t.
$$

As we pointed out above, assuming (\ref{dlil}) uniformly in $\nu
\in {\mathbb Z}$, $f(n_kx)$  satisfies the CLT and replacing
$o(N)$ by $O(N/(\log N)^{1+\varepsilon})$, the LIL also holds for
$f(n_kx)$.  Hence under these conditions  the behavior of
$f(n_kx)$ follows precisely that of i.i.d. random variables.
However, as Fukuyama \cite{fu2} observed, the validity of the CLT
and LIL can break down after a permutation of the terms of the
sequence $(n_k)$, even though an i.i.d.\ sequence remains i.i.d.\
after any permutation. We investigated this surprising phenomenon
in a series of papers \cite{aiberti1,aiberti2,aiberti3}, and found necessary and sufficient Diophantine
conditions for the permutation-invariant behavior of $f(n_kx)$.

In conclusion we note that most results discussed in this chapter
break down for sublacunary sequences (i.e.\ sequences $(n_k)$
satisfying $n_{k+1}/n_k\to 1$), except that the upper half of the
LIL for $f(n_kx)$, i.e.
$$ \limsup_{N\to\infty} \left|\frac{\sum_{k=1}^N f(n_kx)}
{\sqrt{N\log\log N}}\right|<\infty \qquad \text{a.e.}$$ still
holds for some classes of sub-lacunary sequences satisfying strong
number-theoretic conditions (Philipp \cite{pe}, Berkes, Philipp
and Tichy \cite{bpte},  Fukuyama and Nakata \cite{fna},
Aistleitner \cite{adio}; cf. also Furstenberg \cite{furst}, who
studied denseness properties of such sequences from an ergodic
point of view).

The case of superlacunary sequences will be investigated in the
next chapter.

\section{The case $f\not\in L^2$}

In the previous chapter we  saw that under (\ref{f}) and the
Hadamard gap condition (\ref{had}) the asymptotic properties of
partial sums $\sum_{k=1}^N f(n_kx)$  are determined by a
combination of probabilistic and number theoretic effects. In
particular, the behavior of the system $f(n_kx)$ is strongly
influenced by the number of solutions of the Diophantine equation
(\ref{diophequ}). Assuming
\begin{equation}\label{infgap}
n_{k+1}/n_k\to\infty,
\end{equation}
both Diophantine conditions (\ref{dlil}) and (\ref{dlil2}) are
satisfied and consequently $f(n_kx)$ satisfies the central limit
theorem and the law of the iterated logarithm in their classical
form, a result established by Takahashi \cite{tak1,tak2}.
In other words, under the gap condition (\ref{infgap}) the
sequence $f(n_kx)$ behaves precisely as an i.i.d.\ sequence,
without any number theoretic assumptions on $(n_k)$. It is natural
to ask about the asymptotic properties of lacunary sequences
$f(n_kx)$ when the square integrability condition $\int_0^1 f^2(x)
dx<\infty$ does not hold.
Gaposhkin \cite{gapo1968} was the first one to investigate this
question; he proved the following result.

\begin{theorem}\label{gapth}
Let $f: {\mathbb R}\to {\mathbb R}$ be a measurable function with
period 1. Then there exists an increasing sequence $(n_k)$ of
positive integers and measurable functions $g_k(x)$, $\psi_k(x)$,
$\eta_k(x)$, $k=1, 2, \ldots$  on $(0, 1)$ such that the $g_k$ are
stochastically independent and \begin{equation} \label{3terms}
f(n_kx)=g_k(x)+\psi_k(x)+\eta_k(x)
\end{equation} where
\begin{equation}\label{split}
\sum_{k=1}^\infty \|\psi_k\|_M<\infty \quad \text{and} \quad
\sum_{k=1}^\infty \mu \{x: \eta_k (x)\ne 0\}<\infty.
\end{equation}
Here $\| \cdot \|_M$ denotes the norm in the space $M(0, 1)$ of
measurable functions on $(0, 1)$ defined by
$\|\psi\|_M=\inf\{\epsilon>0: \mu (x: |\psi(x)|\ge \epsilon)\le
\epsilon\}$. If $f\in L^p(0, 1)$ $(p\ge 1)$ or $f\in C(0, 1)$,
then the conclusion remains valid with the $g_k$ belonging to the
corresponding spaces and $\| \cdot \|_M$ replaced by  $\| \cdot
\|_p$ or $\| \cdot \|_C$, respectively.
\end{theorem}

 As an immediate consequence, we get
\begin{equation}\label{sumf}
\sum_{k=1}^\infty |f(n_kx)-g_k(x)|<\infty \qquad \text{a.e.}
\end{equation}
in all cases covered by the theorem. Relation (\ref{sumf}) has
powerful consequences. Most limit theorems of probability theory
are invariant for small perturbations, i.e.\ if they are valid for
some sequence $(\xi_k)$ of random variables, then they remain
valid for all sequences $(\xi_k')$ satisfying
$$ \sum_{k=1}^\infty |\xi_k-\xi_k'|<\infty \ \, \text{a.s.}$$
Thus going beyond the limit theorems studied in the previous
section, Ga\-poshkin's theorem extends a very large class of limit
theorems of independent r.v.'s for lacunary sequences $f(n_kx)$.
On the other hand, his theorem provides (except an unproved remark
in the case $f\in L^p(0, 1)$, $p\ge 1$) no estimate of the growth
rate of $(n_k)$ in (\ref{sumf})  and in particular, it provides no
explicit lacunarity condition for many interesting limit theorems
such as versions of the central limit theorem with stable limits,
laws of large numbers and their generalizations, extremal limit
theorems, etc. The purpose of the present chapter is to obtain
explicit growth rates in Gaposhkin's theorem, leading to concrete
lacunarity conditions for a large class of limit theorems for
$f(n_kx)$. As we will see, for many important limit theorems
(including the ones mentioned before) these lacunarity conditions
are actually quite close to (\ref{infgap}). Our results will be
deduced from the following approximation theorem.

\begin{theorem}\label{main} Let $(n_k)$ be an increasing
sequence of positive integers. Then there exists a probability
space $(\Omega, {\cal F}, P)$ and two sequences $(X_k)$ and
$(Y_k)$ of random variables with the following properties.

\medskip\noindent
(a) The sequence $(X_k)_{k\ge 1}$ is a probabilistic replica of
the sequence $\{n_kx)\}_{k\ge 1}$ in the sense that the
distribution of the two sequences in the spaces $(\Omega, {\cal
F}, P)$ and $((0, 1), {\mathcal B}, \mu)$ are the same.

\smallskip\noindent
(b) $(Y_k)_{k\ge 1}$ is an i.i.d.\ sequence with uniform
distribution over $[0, 1]$.

\smallskip\noindent
(c) We have \begin{equation}\label{almiid}
 P(|X_k-Y_k|\ge \delta_k)\le \delta_k  \qquad k=1, 2, \ldots
 \end{equation}
where $\delta_1=1$ and
\begin{equation}\label{delta}
\delta_k=5\left(n_{k-1}/n_k + n_k/n_{k+1}\right) \qquad k=2, 3,
\ldots \, .\\
\end{equation}
\end{theorem}

By the identical distribution of the sequences $(X_k)_{k\ge 1}$
and $(f(n_kx))_{k\ge 1}$ the asymptotic properties of the two
sequences are the same. Using standard probabilistic language, we
can say that we ``redefined'' the sequence $\{n_kx\}_{k\ge 1}$
(without changing its distribution) on the probability space
$(\Omega, {\cal F}, P)$ together with an i.i.d.\ uniform sequence
$(Y_k)_{k\ge 1}$ such that relation (\ref{almiid}) holds with the
$\delta_k$ in (\ref{delta}). In other words, the sequence
$\{n_kx\}_{k\ge 1}$ is, after a suitable redefinition,  a small
perturbation of an i.i.d.\ uniform sequence. This fact implies
that if $\sum_{k=1}^\infty \delta_k<\infty$ (or equivalently if
$\sum_{k=1}^\infty n_k/n_{k+1}<\infty$), then most limit theorems
valid for the i.i.d.\ sequence $(Y_k)_{k\ge 1}$ will be valid for
$\{n_kx\}_{k\ge 1} $ as well. Related, weaker approximation
theorems were obtained in Hawkes \cite{ha} and in Berkes
\cite{b85} dealing with the trigonometric case.

Note that Theorem \ref{main} concerns the specific sequence
$\{n_kx\}$, but depending on the properties of the function $f$,
it leads automatically to a corresponding approximation theorem
for general sequences $f(n_kx)$. For example, if $f$ is continuous
with continuity modulus $\omega(f, \delta)$, then (\ref{almiid})
and (\ref{delta}) imply
\begin{equation}\label{fff}
\sum_{k=1}^\infty |f(X_k)-f(Y_k)|<\infty \qquad \text{a.e.}
\end{equation}
provided
$$
\sum_{k=1}^\infty n_k/n_{k+1}<\infty, \qquad \sum_{k=1}^\infty
\omega(f, n_k/n_{k+1}) <\infty.
$$
A much more general consequence of Theorem \ref{main} is the
following

\begin{theorem}\label{cor0} Let $f: {\mathbb R}\to{\mathbb R}$ be a
measurable function with period 1 and let $(n_k)$ be an increasing
sequence of positive integers. Let $(T_k)$ be positive numbers such
that $\mu\{ x\in (0, 1): |f(x)|\ge T_k\}\le k^{-2}$ and assume
that
\begin{equation}\label{maingapcond}
\sum_{k=1}^\infty \left(T_k \delta_k^{1/4}+\omega_2^{1/2}
(f_{T_k}, 8\delta_k^{1/2})\right)<\infty.
\end{equation}
Then on a suitable probability space there exists a probabilistic
replica $(X_k)_{k\ge 1}$ of $(f(n_kx))_{k\ge 1}$ together with an
i.i.d.\ sequence $(Y_k)_{k\ge 1}$ such that $Y_k$ are distributed
as $f(x)$ on $((0, 1), {\mathcal B}, \mu)$ and
$$\sum_{k=1}^\infty |X_k-Y_k|<\infty \qquad \text{a.s.} $$
\end{theorem}

Here $\delta_k$ is defined by (\ref{delta}),
$$ \omega_2 (f, \delta)= \left(\sup_{0\le h \le \delta} \int_0^1
|f(x+h)-f(x-h)|^2~dx\right)^{1/2}
$$
is the $L^2$ modulus of continuity of $f$ and $f_T$ is the
truncated function $f \cdot \mathds{1} \{|f|\le T\}$. There exist several
classical limit theorems for $f(n_k x)$ involving this modulus of
continuity $\omega_2(f,h)$. For example, Ibragimov \cite{ibra}
proved the CLT for $f(2^k x)$ under (\ref{f}) and the assumption
$$
\sum_{k=1}^\infty \omega_2(f,2^{-k}) < \infty.
$$
Takahashi \cite{tak1} proved the CLT for $f(n_k x)$ under
(\ref{f}), $n_{k+1}/n_k \to \infty$ and
$$
\omega_2(f,h) = \mathcal{O} \left( \log ~\frac{1}{h}  \right)^{-
\alpha}
$$
for some $\alpha > 1$ and Matsuyama and Takahashi \cite{mats}
proved the corresponding LIL under similar, slightly stronger
assumptions. Gaposhkin \cite{gapo1966,gapo1967} proved
that under (\ref{f}) and
\begin{equation} \label{gaposh}
\sum_{k=1}^\infty \omega_2^2(f,n_k/n_{k+1}) < \infty
\end{equation}
the sum $\sum_{k=1}^\infty c_k f(n_k x)$ is a.e. convergent
provided  $\sum_{k=1}^\infty c_k^2 < \infty$ and also that $f(n_k
x)$ satisfies the LIL, provided (\ref{gaposh}) holds with
$\omega_2$ replaced by the ordinary modulus of continuity
$\omega$.

Given any periodic measurable function $f$, we can choose $T_k$ so
that $\mu\{|f|\ge T_k\}\le k^{-2}$ for all $k\ge 1$ and then
condition (\ref{maingapcond}) is satisfied if $\delta_k$ tends to
0 sufficiently rapidly or, equivalently, if $(n_k)$ grows
sufficiently rapidly. More importantly, however, Theorem
\ref{cor0} enables one to give a concrete gap condition implying
the validity of i.i.d.\ limit theorems for lacunary sequences
$f(n_k x)$. We illustrate the procedure on two classical limit
theorems for i.i.d.\ random variables.

\begin{corollary}\label{cor1} Let $f: {\mathbb R}\to{\mathbb R}$ be
a measurable function with period 1 such that the distribution
function
\begin{equation}\label{df} F(x)=\mu\{t \in (0,1): f(t)\le
x\}\end{equation} of $f$ satisfies
\begin{equation}\label{L}
1-F(x)\sim px^{-\alpha} L(x), \qquad F(-x)\sim qx^{-\alpha} L(x)
\qquad \text{as} \ x\to\infty
\end{equation}
for some constants $p, q\ge 0$, $p+q=1$, $0<\alpha<2$ and a slowly
varying function $L$. Let $(n_k)$ be an increasing sequence of
positive integers satisfying (\ref{maingapcond}). Then letting
$S_n=\sum_{k=1}^n f(n_kx)$ we have
\begin{equation}\label{stlim}
(S_n-a_n)/b_n \overset{d}{\longrightarrow} G
\end{equation}
for some numerical sequences $(a_n), (b_n)$ and an $\alpha$-stable
distribution $G$.
\end{corollary}

\begin{corollary}\label{cor2}  Let $f: {\mathbb R}\to{\mathbb R}$
be a measurable function with period 1 such that the distribution
function $F$ in (\ref{df}) satisfies $F(x)<1$ for all $x$ and
$1-F$ is regularly varying at $+\infty$ with a negative exponent.
Let $(n_k)$ be an increasing sequence of positive integers
satisfying (\ref{maingapcond}). Then letting $M_n=\max_{1\le k\le
n} f(n_kx)$, we have
\begin{equation}\label{mn}
(M_n-a_n)/b_n\overset{d}{\longrightarrow}G
\end{equation}
where $G(x)=\exp (-x^{-\alpha}) \mathds{1}_{(0, \infty)}(x)$.
\end{corollary}

Note that (\ref{L}) is the classical necessary and sufficient
condition for the partial sums $S_n$ of an i.i.d.\ sequence with
distribution function $F$ to satisfy the limit theorem
(\ref{stlim})  with suitable norming and centering sequences
$a_n$,  $b_n$. Corollary \ref{cor1} shows that if the distribution
function $F$ of the periodic function $f$ satisfies (\ref{L}),
then the partial sums of $f(n_kx)$ for any $(n_k)$ satisfying
(\ref{maingapcond}) obey the limit theorem (\ref{stlim}).
Similarly, the assumption on $F$ in Corollary \ref{cor2} is the
well-known necessary and sufficient condition for the centered and
normed maxima of an i.i.d.\ sequence with distribution $F$ to
converge weakly to the distribution $G(x)=\exp
(-x^{-\alpha}) \mathds{1}_{(0, \infty)}(x)$, the so called Fr\'echet
distribution. As we know (see e.g.\ \cite{ga}), the limit
distribution in (\ref{mn}) for any i.i.d.\ sequence can be only
one of the Fr\'echet, Weibull and Gumbel distributions with
respective distribution functions $\exp (-x^{-\alpha})\mathds{1}_{(0,
\infty)}(x)$, its analogue on the negative axis and
$\exp(-e^{-x})$; the analogue of Corollary \ref{cor2} holds for
the other two limiting classes, too.

The growth speed of $(n_k)$ in (\ref{maingapcond}) depends on $f$;
clearly, ``nice'' functions $f$ require less rapidly growing
$(n_k)$. For example, in Corollary \ref{cor1} condition (\ref{L})
with $L(x)=1$ and $p=q=1/2$ can be realized  with the function
\begin{equation}\label{repr}
f(x)=
\begin{cases}
-|x-1/2|^{-1/\alpha} \quad &\text{if} \ \ 0<x< 1/2\\
|x-1/2|^{-1/\alpha} \quad &\text{if} \ \ 1/2<x<1.
\end{cases}
\end{equation}
Then we can choose $T_k=k^{1/\alpha}$ and a bound for $|f'|$ on
the set $\{|f|\le T_k\}$ is $Ck^{2(1+\alpha)/\alpha}$ and thus
$$
\omega_2 (f_{T_k}, \delta)\le \omega(f_{T_k}, \delta)\le
Ck^{2(1+\alpha)/\alpha} \delta.
$$
Thus a simple calculation shows that (\ref{maingapcond}) is
satisfied  if
\begin{equation} \label{lac2}
n_{k+1}/n_k \ge k^\gamma
\end{equation}
for some $\gamma=\gamma(\alpha)>0$. Note that (\ref{lac2}) is only
slightly stronger than (\ref{infgap}): relation (\ref{infgap})
requires that $n_k$ grows faster than exponential, while
(\ref{lac2}) is satisfied for $n_k \sim e^{Ck\log k}$ for a
sufficiently large $C$. There are, of course, many other choices
of the function $f$ leading to the same distribution function $F$
in (\ref{df}), which lead in general to faster growing $(n_k)$.

Besides covering a large class of limit theorems, Theorem
\ref{main} leads also to per\-mu\-ta\-tion-invariant results. As
we have seen, under (\ref{f}) and suitable smoothness conditions,
$f(2^kx)$ satisfies the central limit theorem and the law of the
iterated logarithm, but, as Fukuyama \cite{fu2} showed, these
results are not permutation-invariant: both the CLT and LIL break
down after a suitable permutation of the terms of the sequence
$f(2^kx)$. In contrast, relation (\ref{sumf}) is clearly
permutation-invariant and so are its consequences discussed above.

\section{Proofs}

Theorem \ref{main} will be proved by using the strong
approximation technique developed in Berkes and Philipp
\cite{bp1}. More precisely, the result will be deduced from the
following extension of Theorem 2 in \cite{bp1}.

\begin{lemma}\label{proh} Let $(\Omega, {\cal F}, P)$ be a
probability space, ${\cal F}_1 \subset {\cal F}_2 \subset \ldots
\subset {\cal F}$ a sequence of $\sigma$-fields and $X_1, X_2,
\ldots $ a sequence of discrete random variables such that $X_k$
is ${\cal F}_k$-measurable and
\begin{equation} \label{bpp}
 P(\pi \left({\rm dist} (X_k | {\cal F}_{k-1}), {\rm dist}
 (X_k)\right)\ge \gamma_k)
 \le \gamma_k
 \end{equation}
for some $\gamma_k$, $k=1, 2, \ldots$. Assume that on $(\Omega, {\cal F}, P)$ there
exists a random variable $Z$, independent of $\sigma\{ {\cal F}_1,
{\cal F}_2, \ldots \}$ and uniformly distributed over $[0, 1]$.
Then on $(\Omega, {\cal F}, P)$ there exist independent random
variables $Y_1, Y_2, \ldots$ such that $X_k\overset{d}{=}Y _k$
$(k=1, 2, \ldots)$ and
$$P(|X_k-Y_k| \ge \gamma_k)\le 2\gamma_k \qquad (k=1, 2, \ldots).$$
\end{lemma}

Here $\overset{d}=$ denotes equality in distribution, dist$(X_k)$
and dist$(X_k | {\cal F}_{k-1})$ denote, respectively,  the
distribution of $X_k$ and its conditional distribution relative to
${\cal F}_{k-1}$, and $\pi (P_1, P_2)$ denotes the Prohorov
distance of the probability measures $P_1$ and $P_2$ defined by
$$
\aligned \pi(P_1, P_2) = \text{\rm inf} &\bigl\{ \ve > 0 : P_1(A)
\leq
P_2(A^\ve) + \ve \ \text{ and}\\
&\quad P_2(A) \leq P_1(A^\ve) + \ve \text{ for all Borel sets } A
\subset {\mathbb R}\bigr\}
\endaligned
$$
where $A^\ve$ is the open $\ve$-neighborhood of $A$, i.e.,
$$
A^\ve = \bigl\{ x \in {\mathbb R} : |x - y| < \ve \ \text{for
some} \ y \in A\bigr\}.
$$

\bigskip\noindent
{\bf Proof of Lemma \ref{proh}.} We start with recalling some well
known facts from probability theory. Given an atomless probability
space $(\Omega, \mathcal F, P)$ and a distribution function $F$,
there always exists  on $(\Omega, \mathcal F, P)$ a r.v.\ $X$ with
distribution~$F$. As an immediate consequence, if $X$ is a
discrete r.v.\ on an atomless space $(\Omega, \mathcal F, P)$ with
distribution function $F$ and $G$ is a two-dimensional
distribution function with first marginal $F$ (i.e.\ $G(x,
+\infty) = F(x)$), then on $(\Omega, \mathcal F, P)$ there exists
a r.v.\ $Y$ such that the distribution of the vector $(X, Y)$
is~$G$.

Trivially, if $X$ and $Y$ are r.v.'s defined on the same
probability space satisfying $P(|X - Y| \geq \ve) \leq \ve$, then
the Prohorov distance of the distribution of $X$ and $Y$ is $\leq
\ve$. By a theorem of Strassen \cite{st}, the converse is also
true: if $P_1$ and $P_2$ are probability measures on the real line
with $\pi(P_1, P_2) \leq \ve$, then on any atomless probability
space $(\Omega, \mathcal F, P)$ there exist r.v.'s $X$ and $Y$
with distributions $P_1$ and $P_2$ such that $P(|X - Y| \geq \ve)
\leq \ve$. Combining this with the previous remarks it follows
that if the distribution $P_1$ is discrete and $\pi (P_1, P_2)\le
\varepsilon$, then one can even prescribe a r.v.\ $X$ on $(\Omega,
\mathcal F, P)$ with distribution~$P_1$ and there still exists a
r.v.\ $Y$ on $(\Omega, \mathcal F, P)$ with distribution $P_2$
such that $P(|X - Y| \geq \ve) \leq \ve$.

Turning to the proof of Lemma \ref{proh}, we will construct the
r.v.'s $Y_1, Y_2, \ldots$ by induction. We first enlarge the
$\sigma$-fields ${\cal F}_k$ by setting ${\cal F}_k^*=\sigma
\{{\cal F}_k, Z\}$, where $Z$ is
the random variable in the formulation of the lemma. Clearly
dist$(X_k|{\cal F}_{k-1}^*)=$ dist$(X_k|{\cal F}_{k-1})$ and thus
(\ref{bpp}) implies
\begin{equation} \label{bp2}
P(\pi \left({\rm dist} (X_k | {\cal F}_{k-1}^*), {\rm dist}
(X_k)\right)\ge \gamma_k) \le \gamma_k \qquad (k=1, 2, \ldots).
 \end{equation}
Let $Y_1 = X_1$ and assume that $Y_1, Y_2, \dots, Y_{k - 1}$ are
already constructed and satisfy the statements of the lemma,
moreover, $Y_j$ is ${\cal F}_j^*$ measurable for $1\le j \le k-1$.
Since $X_j\overset{d}{=}Y _j$ for $1\le j\le k-1$, the r.v.'s
$Y_1, Y_2, \dots, Y_{k - 1}$ are discrete. Letting ${\cal
F}_{k-1}^{**}=\sigma \{Y_1, \ldots, Y_{k-1}\}$, clearly ${\cal
F}_{k-1}^{**} \subset {\cal F}_{k-1}^*$ and thus (\ref{bp2})
implies
\begin{equation} \label{bp3a}
P(\pi \left({\rm dist} (X_k | {\cal F}_{k-1}^{**}), {\rm dist}
(X_k)\right)\ge \gamma_k) \le \gamma_k \qquad (k=1, 2, \ldots).
\end{equation}
Consider the sets $D$ of the form
\begin{equation}\label{D}
D = \left\{ Y_1 = b_1, \dots, Y_{k - 1} = b_{k - 1} \right\}
\end{equation}
where $b_1, \ldots, b_{k-1}$ are in the range of $Y_1, \ldots,
Y_{k-1}$, respectively. The union of these sets is clearly
$\Omega$; we will construct $Y_k$ on each such set separately. We
clearly have  dist$(X_k|{\cal F}_{k-1}^{**})=$ dist$(X_k|D)$ on
the set $D$ in (\ref{D}), and thus (\ref{bp3a}) implies that the
sets $D$ can be distributed into two classes $\Gamma_1$ and
$\Gamma_2$ such that $\sum_{D\in \Gamma_2} P(D)\le \gamma_k$ and
for any $D\in \Gamma_1$ we have
\begin{equation}\label{pi}
 \pi (\hbox{dist} (X_k|D), \hbox{dist}(X_k))\le
\gamma_k.
\end{equation}
Let first $D\in \Gamma_1$, let $P^{(D)}$ denote conditional
probability with respect to $D$  and define the probability
measures $P_1$ and $P_2$ on the Borel sets of $\mathbb R$ by
$$
P_1 (A)=P(X_k\in A|D), \quad P_2 (A)=P(X_k\in A),\qquad  A\in
{\cal B}.
$$
By (\ref{pi}) we have $\pi (P_1, P_2)\le \gamma_k$ and since the
probability space $(D, {\cal F}_k^*, P^{(D)})$ is atomless,
the remarks at the beginning of the proof imply that on the
probability space $(D, {\cal F}_k^*, P^{(D)})$ there exists a
random variable $Y_k$ with distribution dist$(X_k)$ such that
$P^{(D)} (|X_k-Y_k|\ge \gamma_k)\le \gamma_k$. This defines the
random variable on each set $D\in \Gamma_1$; for a set $D\in
\Gamma_2$ let $Y_k$ be any r.v.\ on the probability space $(D,
{\cal F}_k^*, P^{(D)})$ with distribution dist$(X_k)$. Thus we
defined $Y_k$ on the whole probability space $\Omega$; clearly,
the so defined $Y_k$ is ${\cal F}_k^*$ measurable and has the
property that its conditional distribution on any set of the form
(\ref{D}) equals its unconditional distribution and thus $Y_k$ is
independent of the vector $(Y_1, \ldots, Y_{k-1})$. Also,
$$
P(|X_k-Y_k|\ge \gamma_k)\le 2\gamma_k.
$$
This completes the induction step and thus the proof of Lemma
\ref{proh}.

\bigskip\noindent {\bf Proof of Theorem \ref{main}.}
Let $\ve_k = n_k/n_{k + 1}$ and $\psi(x)= \{x\}$. Set $A_k =
\bigl\{i/n_{k + 1} : 1 \leq i \leq n_{k + 1}\bigr\}$, $B_k =
\bigcup^k_{j = 0} A_j$ and let $J_1, J_2, \dots, J_{m_k}$ be the
left closed intervals to which the points of $B_k$ divide the
interval~$[0, 1)$. Let $\mathcal F_k$ be the $\sigma$-field
generated by the intervals $J_1, \dots, J_{m_k}$ and set
$$
\aligned
T_k &= T_k(x) = \psi(n_k x),\\
X_k &= E(T_k | \mathcal F_k).
\endaligned
$$
Clearly $\mathcal F_1 \subset \mathcal F_2 \subset  \cdots$ and
$X_k$ is an $\mathcal F_k$-measurable discrete r.v.\ taking the
constant value $\mu (A)^{-1} \int_A \psi(n_k x)\, dx$ on any atom
$A$ of~$\mathcal F_k$. Since the atoms of $\mathcal F_k$ are
intervals of length $\leq 1/n_{k + 1}$, we have
\begin{equation}\label{eq:13}
|T_k - X_k| \leq n_k / n_{k + 1} = \ve_k\, .
\end{equation}
Let $\mathcal G_k$ denote the set of intervals $I = \bigl[ i /
n_{k + 1},  (i + 1)/n_{k + 1}\bigr),~0 \leq i \leq n_{k + 1} - 1$
which contain in their interior no point of $B_{k - 1}$. Clearly
each interval $I \in \mathcal G_k$ is an atom of $\mathcal F_k$
and
\begin{equation}
\sum_{I \in {\mathcal G}_k} \mu (I) \geq 1 - 2\ve_k \quad \ k \geq
k_0. \label{eq:144}
\end{equation}
To see the last inequality observe that the number of those
intervals
$$I = \bigl[i/n_{k + 1}, (i + 1) / n_{k + 1} \bigr),~0
\leq i \leq n_{k + 1} - 1,
$$
which contain in their interior a point of $B_{k - 1}$ is at most
$\text{\rm card}\, B_{k - 1} \leq n_1 + \dots + n_k$ and thus the
total measure of these intervals is at most $(n_1 + \dots + n_k) /
n_{k + 1}$. Now (\ref{infgap}) implies $n_1 + \dots + n_k \leq
2n_k$ for $k \geq k_0$, proving \eqref{eq:144}.

Let $H$ denote the uniform distribution over $(0,1)$. Clearly the
conditional distribution of $T_k = \psi (n_k x)$ relative to each
interval $I \in {\mathcal G}_{k-1}$ is $H$ and thus by
\eqref{eq:13} we have
\begin{equation*}
\pi\left( \text{dist} (X_k|I), H\right) = \pi\left( \text{dist}
(X_k|I), \text{dist}(T_k|I)\right)\le \e_k \qquad \text{for} \
I\in {\mathcal G}_{k-1}.
\end{equation*}
In view of (\ref{eq:144}), the last relation means
\begin{equation} \label{pr2}
\pi\left( \text{dist} (X_k|{\cal F}_{k-1}), H\right) \le \e_k
\quad \text{with probability}  \ \ge 1-2\ve_{k-1}.
\end{equation}
Since dist$(X_k)$ is obtained from dist$(X_k|{\cal F}_{k-1})$ by
integration, (\ref{pr2}) easily implies
\begin{equation*}
\pi\left( \text{dist} (X_k), H\right) \le \e_k + 2\ve_{k-1},
\end{equation*}
which, together with (\ref{pr2}), yields
\begin{equation*}
\pi\left( \text{dist} (X_k|{\cal F}_{k-1}),
\text{dist}(X_k)\right) \le 2(\e_k+\e_{k-1}) \quad \text{with
probability} \ \ge 1-2(\ve_k+ \ve_{k-1}).
\end{equation*}
We thus showed that condition (\ref{bpp}) of Lemma \ref{proh}
holds with $\delta_k = 2(\ve_k + \ve_{k-1})$. Hence if on the
probability space $((0, 1), {\cal B}, \mu)$ there exists a
uniformly distributed random variable $Z$ independent of $\sigma
\{ {\cal F}_1, {\cal F}_2, \ldots, \}$, then
Lemma \ref{proh} applies and yields a sequence $(Y_k)$ of
independent random variables on this space such that $Y_k
\overset{d}{=} X_k$ and
$$
\mu \bigl(|X_k - Y_k| \geq \delta_k\bigr) \leq 2\delta_k \qquad \
k \geq k_0.
$$
By (\ref{eq:13}) this implies
$$
\mu \bigl(|\{n_kx\} - Y_k| \geq \delta_k'\bigr) \leq \delta_k'
\qquad \ k\ge k_0
$$
with $\delta_k'=5(\ve_k + \ve_{k-1})$. If a random variable $Z$
with the desired properties does not exist, replace $((0, 1),
{\cal B}, \mu)$ by the product space $((0, 1), {\cal B}, \mu)
\times ((0, 1), {\cal B}, \mu)$ and redefine all random variables
and $\sigma$-fields on the new space in an obvious fashion. In the
new space the required $Z$ obviously exists and all arguments of
our proof remain valid in the new space.

To complete the proof of Theorem \ref{main} it suffices now to
show that there exists an i.i.d.\ sequence $(Z_k)$ with $Z_k
\overset{d}{=} \psi$ such that $|Y_k-Z_k| \leq \ve_k$ for $k=1, 2,
\ldots$.
By passing to a suitable product space as before,
we can assume without loss of generality that there exist
independent random variables $\eta_1, \eta_2, \ldots$, having
uniform distribution over $(0, 1)$ and independent also of $X_1,
Y_1, X_2, Y_2, \ldots$; let ${\cal H}_k = \sigma\{ Y_k, \eta_k\}$.
By the remarks made at the beginning of the proof, on the atomless
probability space $(\Omega, {\cal H}_k, P)$ there exists a random
variable $Z_k$ such that the joint distribution of $Y_k$ and $Z_k$
is the same as the joint distribution of $Y_k$ and $\psi$ in $(
(0, 1), {\cal B}, \mu)$. Clearly the $Z_k$ are independent, $Z_k
\overset{d}{=} \psi$  and by (\ref{eq:13}) we have $|Z_k-Y_k|\le
\ve_k$. This completes the proof of Theorem
\ref{main}.\\

\bigskip\noindent
{\bf Proof of Theorem \ref{cor0}.} We first consider the case when
$f$ satisfies the additional condition  $\int_0^1 f^2(x)\,
dx<\infty$;
let
$$ f\sim \frac{a_0}{2}+\sum_{k=1}^\infty (a_k \cos 2\pi k x+b_k \sin 2\pi kx)
$$
be its Fourier series. Then for any $n\ge 1$ we have
$$ \|f-s_n(f)\|\le \omega_2(f, \pi/n)$$
where $s_n(f)$ denotes the $n$th partial sum of the Fourier series
of $f$ (see e.g.\ \cite{ba}, Vol.\ II, p.\ 156).
Write $f=f_1+f_2$, where
$$f_1(x)=s_m (f, x)=\frac{a_0}{2}+\sum_{k=1}^m (a_k \cos 2\pi k x+b_k \sin 2\pi kx),
$$
and $f_2=f-s_m$ with $m$ to be determined later. By Theorem
\ref{main} there exists a probability space $(\Omega, {\cal F},
P)$ and sequences $(X_k^*)$ and $(Y_k^*)$ of random variables such
that $(X_k^*)$ is a probabilistic replica of the sequence
$\{n_kx\}$, $Y_k^*$ are i.i.d.\ random variables with uniform
distribution over $(0, 1)$ and
\begin{equation}\label{almiid2}
P(|X_k^*-Y_k^*|\ge \delta_k)\le \delta_k  \qquad k=1, 2, \ldots,
 \end{equation}
with $\delta_k$ defined by (\ref{delta}). Then
$$
|f(X_k^*)-f(Y_k^*)|\le |f_1(X_k^*)-f_1(Y_k^*)|+
|f_2(X_k^*)-f_2(Y_k^*)|=: V_1+V_2.
$$
Clearly
\begin{align*}
&|f_1'(x)| \le 2\pi \sum_{k=1}^m k(|a_k|+|b_k|)\\
&\le 2\pi \left[\left(\sum_{k=1}^m a_k^2\right)^{1/2} +
\left(\sum_{k=1}^m b_k^2\right)^{1/2}\right] \left(\sum_{k=1}^m
k^2\right)^{1/2}\le 4\pi\|f\|m^{3/2}
\end{align*}
and thus (\ref{almiid2}) and the mean value theorem imply that
$|V_1|\le 4\pi\|f\|m^{3/2} \delta_k $ with probability $\ge
1-\delta_k$. On the other hand, $X_k^*$ and $Y_k^*$ are uniformly
distributed r.v.'s over $(0, 1)$ and thus
$\|f_2(X_k^*)\|=\|f_2(Y_k^*)\|=\|f_2\|$.  Hence we have
$$\|V_2\| \le 2 \|f_2\|\le 2\omega_2 (f, \pi/m)$$
and thus the Markov inequality yields that $|V_2|\le
\omega_2^{1/2} (f, \pi/m)$ with the exception  of a set with
measure not exceeding $4\omega_2 (f, \pi/m).$  We thus proved that
\begin{equation}\label{ut}
|f(X_k^*)-f(Y_k^*)|\le 4\pi\|f\|m^{3/2} \delta_k + \omega_2^{1/2}
(f, \pi/m)
\end{equation}
with probability exceeding $1-(\delta_k+ 4\omega_2 (f, \pi/m))$.
In the case of a general periodic measurable $f$, apply the
previous argument for the truncated function $f_{T_k}=f \cdot \mathds{1}\{|f|\le
T_k\}$, with the choice $m=[\delta_k^{-1/2}]$. Clearly,
$\|f_{T_k}\|\le T_k$, and thus using (\ref{ut}) for $f_{T_k}$ and
applying the Borel-Cantelli lemma, Theorem \ref{cor0} follows with
$X_k=f(X_k^*)$, $Y_k=f(Y_k^*)$.


\end{document}